\documentclass{article}
\usepackage[margin=3cm]{geometry}
\usepackage{amsfonts,amssymb,amsmath,amscd}
\usepackage{titlesec}
\usepackage{booktabs}
\usepackage{array}
\usepackage{longtable}
\usepackage{hyperref}
\titleformat*{\section}{\large\bfseries}

\title{Computational Platonism}

\author{Asvin G\thanks{This project was supported by a grant from the Simons Foundation International [MPS-SIM-00001691, JT]. The author thanks Claude Opus 4.7 and Fable 5.0 (Anthropic) for substantial editorial collaboration on this essay.}\\
Institute for Advanced Study, Princeton}

\begin{document}
\maketitle

\begin{quote}
\centering
 \textit{The product of mathematics is clarity and understanding. Not theorems, by themselves.}

 \raggedleft{Thurston, Proof and Progress}
\end{quote}

Many mathematicians will agree with Thurston without being able to say \textit{what} they are understanding. Other sciences have a somewhat easier answer: physicists explain matter, biologists study life, astronomers describe the night sky; each points to a domain of phenomena that is observable, reproducible, and independent of any particular theory about it. What can mathematicians point to? Davis and Hersh \cite{davis2012mathematical} proposed ``mental objects with reproducible properties''. Arnold \cite{arnold1997teaching} went further, calling mathematics ``the part of physics where experiments are cheap''. Both catch genuine features of mathematical practice, and both leave something essential out: Davis and Hersh keep mathematics trapped in the mind, while Arnold leaves open the nature of its experiments.

I will argue that mathematics is, at its core, an experimental science.\footnote{The view I develop has antecedents in the quasi-empiricism of Lakatos \cite{lakatos1963proofs} and Chaitin \cite{chaitin2005metamath}, in Putnam's empirical account of axioms \cite{putnam1975mathematical}, in Maddy's naturalism about axiom selection \cite{maddy2011defending}, and in Wolfram's program of treating computation as foundational \cite{wolfram2002}. It differs from the humanism of Hersh \cite{hersh1997what} in grounding mathematical objectivity in substrate-independent computation rather than in mathematical practice.} Its experiments are computations: not arithmetic in the schoolroom sense but computation as Turing formalized it, a class of physical processes, carried out by definite rules, whose outcomes are reproducible and can surprise the person who initiates them. Like any experiment, a computation is something you \textit{do} and its result something you \textit{observe}. Mathematicians observe patterns in these experiments and construct axiomatic frameworks to explain them, as physicists construct theories to explain patterns in physical experiments. And like any scientific theory, an axiomatic framework can be tested, found incomplete, or even contradicted, a perspective that will lead to a new reading of G\"odel's incompleteness theorems.

The best scientific theories go beyond explaining observed phenomena on their own terms. They posit entities that are not directly observable but seem to causally determine what we do observe: Maxwell's electromagnetic field, which unifies Faraday's observations, and Darwin's heritable units (what we now call genetic material), without which natural selection cannot operate. Mathematical theories often play a similar role, and the entities they reach for are just as dramatic: completed infinities, the continuum as a single totality, higher-dimensional geometries. We intuit what our exploration seems to point to, and give each intuition computational flesh by specifying an axiomatic framework.

The Platonism of the title refers to this reaching, and goes with a revisionist reading of Plato that puts the inference center stage over the Forms, as will become clear.

This contrasts with the common view of mathematics as a purely deductive science, in which axioms are the foundation and theorems are derived from them, so that computation is a tool for exploring consequences. Here the relationship is reversed: computation is the foundation, and axioms are the theory we build to explain it. Making this precise requires first saying what a computation is.

\section*{\centering{Computational traces}}

Mathematicians produce tables of primes, lists of factorizations, records of which integers are sums of two squares. They verify proofs step by step and derive theorems from axioms. These are outputs of definite procedures, and we will call them \textit{computational traces}: the raw data of mathematical experience, the analogue of readings from a physical experiment. What makes them experimental? They are reproducible, they can surprise us, and they are durable, outliving the theories that interpret them.

A calculation by one person can be checked by another and the results will agree; a proof can be verified by anyone who follows the steps; the hundredth prime is the same whether computed by hand or by machine. What does this reproducibility require?

Consider the analogous property in a physical experiment. In a scattering experiment, the outcome we record is the pattern on the detector screen, not the full microscopic state of every particle: the exact trajectories, the thermal state of the detector, the time of day are all discarded. This is an idealization, and only relative to it does reproducibility become a meaningful claim. Reproducibility is then a statement about the idealized outcome under an assumed symmetry: the experiment gives the same result in different laboratories and at different times. The symmetries are themselves empirical claims; Galileo's results on free fall hold on the Earth's surface but not in orbit. Without the idealization there is no definite outcome to reproduce; without the symmetry there is no reason to expect it to recur.

Mathematical reproducibility has the same structure. The idealization: we care about the symbolic content of a calculation, not the ink, the handwriting, or the neurons that fired. Davis and Hersh were right to the extent that this idealization makes mathematical objects seem mental. But the process producing the symbols is physical, whether on paper, on a circuit board, or in the mind; even mental arithmetic implements a procedure a machine could equally carry out. The symmetry: the procedure gives the same result regardless of substrate, and this invariance grounds mathematical reproducibility in the physical world. That these symmetries are so much harder to imagine failing than their physical counterparts is part of what makes mathematical knowledge seem so certain.

Turing's formalization of computation \cite{turing1936computable} made both precise. His machine reads symbols on a tape one at a time and writes new ones according to a finite table of rules; it is specified entirely by that table and the initial tape. Despite this simplicity the setup is universal: any procedure that follows definite rules, on any physical system whatsoever, can be carried out by such a machine. The idealization is now explicit, the outcome being what is on the tape when the machine halts, and so is the symmetry: any physical system implementing the same rules on the same input produces the same output.

A mathematical \textit{definition} specifies a computational system: a string and a set of rewrite rules operating on it. The natural numbers, for example, are specified by an initial symbol $\underline{0}$ and operations $S,+,\times$. The $n$-th natural number is the outcome of applying $S$ to $\underline{0}$ $n$ times, $\underline{n} = S^n\underline{0}$, and $+,\times$ are defined by rewrite rules:
\begin{enumerate}
    \item $\underline{n} + \underline{0} \to \underline{n}$
    \item $\underline{n} + S\underline{m} \to S(\underline{n}+\underline{m})$
    \item $\underline{n} \times \underline{0} \to \underline{0}$
    \item $\underline{n} \times S\underline{m} \to \underline{n} + \underline{n}\times \underline{m}$.
\end{enumerate}
Given any two natural numbers as input, these rules produce a definite output. They are a procedure, not a claim.

A very different example is first-order logic. The strings are formulas built from variables, connectives, and quantifiers; the rewrite rules are the rules of deduction: from ``$P$ implies $Q$'' together with $P$, derive $Q$. A proof is a sequence of such rewrites from an initial set of formulas called axioms. Like arithmetic, this is a computational system of definite rules on strings, but its outputs are theorems rather than numbers.

Computational traces are not merely reproducible; they can surprise us. The rewrite rules of arithmetic generate, among other things, the perfect squares $1, 4, 9, 16, 25, \ldots$ The number $\pi$, the ratio of a circle's circumference to its diameter, is a geometric quantity with no apparent connection to whole numbers. Yet Euler showed in 1734 that $1/1 + 1/4 + 1/9 + 1/16 + \cdots$ equals exactly $\pi^2/6$: geometry emerged from an arithmetic computation. Gauss observed in the 1790s that the proportion of primes near a large integer $X$ seemed to be roughly $1/\log X$, a statistical pattern in prime tables that was compelling long before any theory existed to explain it.

In both cases the global behavior of the system outruns what its rules make obvious. Each step is routine, but the patterns across many steps are not foreseeable without running the system. Computation is experimental in exactly this sense: you have to do it to find out what happens.

The same happens throughout mathematics. Euler counted the vertices, edges and faces of many polyhedra and observed that $V - E + F = 2$, a topological invariant emerging from simple counting, and the central case Lakatos used in \textit{Proofs and Refutations} \cite{lakatos1963proofs} to develop his quasi-empirical view. Fourier assumed in 1822 that every function could be expressed as a trigonometric series and found that this unjustified assumption produced brilliant, verifiable solutions of the heat equation, sparking the questions that led to Cantor's set theory and the foundational crisis. In each case the computational system was richer than anyone anticipated from its defining rules alone.

\begin{quote}
\centering
 \textit{Mais les objets sont têtus et les méthodes explicites ne cessent de ressurgir. Un calcul est toujours plus général que le cadre théorique dans lequel on l'enferme à une période donnée.}

 \smallskip

 \textnormal{But objects are stubborn and explicit methods never cease to resurface. A calculation is always more general than the theoretical framework in which one confines it at any given period.}

 \raggedleft{M. Demazure, \textit{Résultant, Discriminant}}
\end{quote}

Finally, traces are durable in a way their interpretations are not. Physical phenomena are reinterpreted as foundational theories evolve; the nature of light was debated for centuries before quantum electrodynamics. Fermat's principle, that light follows the path of least time, was initially seen as a mathematical curiosity while the particle view was dominant, since read literally it requires either that light be a wave or that its particles have foreknowledge of the future, in violation of causality. The wave theory, from Huygens to Fresnel, reinterpreted it as interference: paths near the shortest reinforce one another while the rest cancel. Feynman reinterpreted it once more as the photon taking all paths. The principle is now fundamental physics, but its interpretation has evolved with our understanding of the underlying system. The same happens in mathematics, perhaps to a greater extent.

The Riemann--Roch theorem is an example. Riemann read it as a statement about the number of linearly independent meromorphic functions on a Riemann surface with specified poles. It was later reread by Brill, Noether, and the Italian school, and made rigorous by Zariski, as a statement about complete linear systems of divisors on algebraic curves and surfaces, and then once more by Serre and Grothendieck as a statement about the Euler characteristic of a coherent sheaf. The original observation is an experimental fact about a pattern in the computational traces of algebraic geometry, and it has been the kernel of each interpretation, pointing towards deeper understanding at each stage.

This may seem close to formalism, which also places symbolic manipulation at the center of mathematics. But for the formalist the formal system is all there is; there is nothing to be surprised \textit{by}. Here the traces a formal system generates are the subject matter of mathematics, not its entirety: they are phenomena that demand explanation. At this pre-axiomatic level mathematics is completely constructive and unambiguous; we learn to add before we have any idea of induction, let alone the Peano axioms. If traces are the experimental data, the natural question is what plays the role of theory.

\section*{\centering{Axioms as theory}}

In physics, a theory predicts the outcomes of experiments not yet performed and explains why past ones came out as they did. Newton's theory predicted the existence and orbit of Neptune and explained why planets move on ellipses. Prediction lets us replace experiment with computation; explanation reveals structure we would not otherwise see. Theories also point to new targets, phenomena nobody was looking for before. Dirac's relativistic equation for the electron had negative-energy solutions that could not be dismissed as artifacts; he predicted these corresponded to a new particle, the positron, discovered shortly after. The equation did not just explain the electron; it revealed antimatter.

For the mathematical parallel, take Fermat's Last Theorem. The conjecture is a pattern in the traces of the natural numbers: for $n > 2$, no computation of $a^n + b^n$ produces a perfect $n$-th power. Small cases can be checked directly, and the first proofs, for $n = 3$ and $n = 4$, work within elementary number theory using the Peano axioms. But already in the 1840s Kummer found that going further required expanding the computational framework. He needed unique factorization in rings of cyclotomic integers, and to his surprise it fails in general. The failure forced him to invent ``ideal numbers,'' a new kind of computational object motivated not by abstract considerations but by the desire to recover, in larger number systems, a property the ordinary integers possess. The definition, later recast by Dedekind as that of an ideal, is far from obvious; it was shaped by what Kummer saw when he computed in these rings.

The expansion continued for well over a century. Each new computational system along the way (algebraic number fields, elliptic curves, modular forms, Galois representations) arose in response to patterns observed in existing ones. Galois found that polynomial equations were better explained by their symmetry groups than by explicit solutions, in the very effort to find such solutions. Grothendieck rebuilt algebraic geometry on sheaf-theoretic foundations because the targets of number theory required new cohomological tools. Wiles's proof of Fermat's Last Theorem draws on this entire apparatus, vast machinery bearing down on a question stated purely in terms of the natural numbers and their rewrite rules. Each new framework is at once a theory explaining patterns in the traces of the previous one and a new computational system generating traces for the next theory to explain.

On the standard view, axioms are foundational truths and need no defense. On the view developed here they are theoretical claims about computational systems, claims that could turn out to be wrong. What then justifies our confidence in Peano arithmetic (PA) or Zermelo--Fraenkel set theory with choice (ZFC)? The same thing that justifies confidence in a physical theory: each has predicted and explained a mass of evidence without producing a contradiction. We believe the Peano axioms because they have never led to a false prediction about the natural numbers, and ZFC because it has organized and explained vast swaths of mathematics without inconsistency. This is not certainty but the provisional confidence physicists place in general relativity, the kind of justification Maddy has developed in detail \cite{maddy2011defending}. Since it rests on intuitions about the underlying computational systems, the next question is what those intuitions reach for.

\section*{\centering{What axioms reach for}}

As frameworks grow more powerful, their fundamental objects tend to move beyond direct observation or computation. The basic quantities of classical mechanics (position, velocity, force) correspond straightforwardly to what we measure. The wavefunction does not; it determines probabilities through the Born rule but is never itself observed. The metric tensor of general relativity describes spacetime geometry but cannot itself be measured, a consequence of general covariance. This move in physics is a strict parallel of the move in mathematics from constructive to non-constructive, and in many ways the bread and butter of the subject.

Arithmetic only ever produces specific natural numbers $\underline{n}$, yet any investigation of it at once suggests concepts beyond it: the totality of natural numbers as a set, modular arithmetic, the real numbers as a continuum. These are intuited objects, posited because the patterns in finite traces look like shadows of something more general, and they earn their place as the wavefunction and the metric tensor do: through what they let us prove and explain about the things we can directly compute with.

The cleanest historical case is the parallel postulate. For over two millennia it was treated as a self-evident axiom of geometry, a claim about the computational system of ruler and compass, much as induction is treated today as a self-evident truth about the system of the successor function. Belief in it was grounded in direct experience of physical geometry, but as Lobachevsky, Bolyai, and Riemann showed in the nineteenth century, geometries that violate it can be conceived, and were in fact found by mathematicians working with its negations and building rich intuitions about the resulting systems. Both the belief and the expansion arose from computational experience of the kind described here. We are living through the same episode with the continuum hypothesis, which the axioms of set theory neither prove nor refute because the intuitions they codify, drawn from finite collections, say nothing about how many subsets an infinite set has.

Even concepts as simple as negative numbers arose this way. Until the mid-sixteenth century, European arithmetic was limited to the positive rationals; numbers measured lengths or counted quantities, and negative numbers were patently ridiculous under this framework. Yet in solving polynomial equations with natural-number coefficients, mathematicians needed intermediate expressions involving negative numbers, and sometimes imaginary numbers like $\sqrt{-1}$, even when the final answer was a natural number. Arithmetic was not rich enough to capture the patterns they were seeing, so they expanded it by positing new entities that could be computed with. The integers, and later the complex numbers, are now so familiar that we forget how they came about, but each is a prototype of the process described here.

The split between constructive and non-constructive mathematics is another instance, though not always recognized as one. Constructive mathematics, while not precisely delimited, can be seen as a peak of empiricism: we demand explicit constructive certificates for any claim about a computational system such as arithmetic. Non-constructive mathematics arises when we make precise our intuitions about entities beyond the original system by building new axiomatic frameworks for them. The infinite cannot be grasped directly within the natural numbers, but we can make it accessible to computation by encoding the properties of it we can intuit. The axiom of choice, the prototypical non-constructive principle, codifies particular intuitions about the behavior of infinite sets, intuitions learned from experience with finite collections but not logically necessary, and the codification can be made precise and explored for consistency and consequences like any other computational system.

This reaching is the sense in which the title speaks of Platonism. In the allegory, prisoners mistake the shadows on the cave wall for the objects casting them, and the philosopher's task is to escape the cave and confront the objects. I keep the inference from shadow to shape but give up the reality at the end of it. Platonism, as I use the word, is not the claim that a realm of perfect objects stands behind the shadows; it is the activity of inferring a shape from its shadows, a compression in which many scattered traces are explained at once by a single posited object. The parallel postulate, the negative numbers, the axiom of choice, and the axiom of infinity were each such an inference.

\section*{\centering{Axioms as tools for exploration}}

What the original reading of the allegory misses is that the inference never terminates. Each inferred shape, once given computational flesh, is a computational system in its own right, casting new shadows from which we infer further shapes.

The axioms specify the rules; the rules generate traces; the traces become the next layer of evidence to explain. Kummer's pursuit of Fermat's Last Theorem ran through this loop several times. Recall that he moved into the cyclotomic integers $\mathbb{Z}[\zeta_p]$, where $a^p + b^p$ splits as $\prod_{k=0}^{p-1} (a + \zeta_p^k b)$, found that unique factorization fails there, and repaired it by positing ideal numbers: abstract objects, not elements of the ring, that behave like factors and restore unique factorization one level up. Once defined, their computational properties could be explored on their own terms. Ideal numbers modulo actual ones form a group, itself far from obvious, and looking for patterns in the size of this class group Kummer came to a surprising $p$-adic connection between these class groups and values of the zeta function at negative integers. He had not set out to find it; it emerged from the structure he had instantiated, and led a century later to the deep work of Iwasawa, Mazur--Wiles, and many others.

A different way axiomatization powers exploration is by unifying seemingly disparate phenomena and suggesting new examples. Take the Weil conjectures. Counting the solutions of a polynomial equation over finite fields produces a sequence of integers, one for each field size, as concrete a trace as arithmetic offers. Reading Gauss's computations of this kind, Weil observed in 1949 that these sequences behave as if governed by a hidden topology: they satisfy a rationality property, a functional equation, and an analogue of the Riemann hypothesis, exactly as they would if varieties over finite fields had a cohomology theory like that of a complex manifold. No such theory existed, and it was hard to imagine how discrete spaces like finite fields could carry topological structure. The conjecture was the inference of a shape from its shadows, and Grothendieck and Serre built the shape: étale cohomology, reached by first recasting cohomology itself in terms of abelian categories and derived functors, a recasting general enough that group cohomology, the cohomology of complex manifolds, and the new arithmetic cohomology became instances of one construction. Deligne proved the last of the conjectures inside it in 1974. The posit did not stop at point counts. The same machinery now organizes Galois representations and modular forms and carried the proof of Fermat's Last Theorem, and Quillen's extension of it to homotopical settings opened the road to modern homotopy theory and the $\infty$-categories beneath it.

This dual power, to explain existing phenomena and to serve as the basis for further exploration, is what makes axiomatic frameworks so central to mathematics. They are usually grounded in traces of familiar objects, which is why they tend to be resilient. Each is a creative, non-algorithmic leap, disciplined by the traces it aims to explain and verified by practice with the resulting system.

The parallel with physics is most interesting precisely here. In physics, theories explain observations and suggest new experiments but never substitute for experiment itself. In mathematics the distinction collapses: the theories are themselves computational systems, and computing in them is what we mean by experiment. The language of explanation and the language of discovery coincide. This is the source of the sharp difference in flavor between the two fields, and the reason the experimental nature of mathematics is so often overlooked. Even so, mathematical theories are as provisional as physical ones, and as subject to revision and replacement when they fail to capture the phenomena.

\section*{\centering{How theories fail}}

Consider Peano arithmetic again. Its axioms play two roles. The rewrite rules for $S, +, \times$ are \textit{definitional}: they specify a procedure that simply runs. The induction axiom is \textit{descriptive}: it claims that every property holding of $\underline{0}$ and closed under $S$ holds of all natural numbers, a claim beyond any finite verification. The $=$ sign marks the boundary: in $\underline{n} + \underline{0} \to \underline{n}$ the arrow defines a procedure; in ``for all $n, m$: $n + m = m + n$'' equality asserts a claim. The distinction determines where contradictions can arise. Rewrite rules cannot contradict themselves: a computation is a physical process and its result is whatever the process produces. Substrate-independence could in principle fail, but that would be a physical violation, not a mathematical contradiction. Contradictions arise in the descriptive claims, such as the principle of induction: assertions about the system that turn out to be incompatible, as Newton's theory turned out to be incompatible with the perihelion of Mercury.

Set theory is a case where the practice being formalized already contained contradictions. Cantor's analysis of the continuum and everyday intuitions about finite sets converged on a few implicit assumptions, chief among them that any property determines the set of things satisfying it. Analysts had used this principle heavily and to good effect, but as Russell showed, applied to the totality of all sets it contradicts itself: the set of all sets that do not contain themselves contains itself exactly when it does not, in barely a line. Naive set theory seemed hopelessly inconsistent.

The task was therefore not to describe a consistent system precisely, as with arithmetic, but to repair an inconsistent practice: to find consistent axioms matching our intuitions about finite collections as closely as possible in the infinite regime. The result of this search was the canonization of ZFC as the backbone of modern mathematics. As far as we can tell it is consistent, and powerful enough to describe essentially all of mathematics, but it is incomplete on questions of natural interest such as the continuum hypothesis. Where PA proves essentially any theorem a number theorist would ask about, ZFC is an incomplete theory matching our incomplete intuitions about the infinite. Set theorists have taken this as an invitation to explore the wild world ZFC describes, uncovering surprising phenomena such as large cardinals: axioms built to explain one practice became the instrument of a new one.

A more striking case of intuition outrunning the computational substrate is early algebraic geometry. The Italian school achieved remarkable results about curves and surfaces through geometric intuition disciplined by examples, but the intuition was never given definite computational form; their arguments relied on what happens ``in general,'' a notion no one had a procedure for. The first generation, Castelnuovo and Enriques, had extraordinary intuition; as Mumford argues in a careful reassessment \cite{mumford2011intuition}, Enriques's ideas about infinitesimal deformations essentially anticipated Grothendieck's schemes decades before the tools existed to make them precise. But the style was fragile. Severi pushed it further and went astray, claiming in 1934 that the space of rational equivalence classes of zero-cycles on a smooth surface depends on only finitely many parameters; in 1968 Mumford showed that for surfaces of a certain type the space is infinite-dimensional. Severi had made no arithmetic error, and his computations with specific surfaces were correct. The system he had in mind was consistent; it just was not the system the computations were tracing out. The response, carried out over decades by Zariski, Weil, and Grothendieck, was not to restrict the framework but to enlarge it, until singularities, families, moduli, and ``general position'' each had a precise procedural meaning.

Since axiomatic frameworks are theories about an underlying physical process (the computational system and its traces), a mistaken intuition encoded in the axioms can be corrected, as in physics. Inconsistencies are not fatal to mathematics in our framework but signposts toward a more accurate explanation of the computational processes we do have access to. Physics faces the same problem and treats each theory as provisional without illusions of finality; even so, the dream of a complete theory of all physical processes remains alive; we are far from such a theory currently, but there is no reason to think it impossible. Mathematics finds itself in a very different place.

\section*{\centering{Incompleteness as inexhaustibility}}

Hilbert's program aimed for axioms complete enough to derive every truth about the natural numbers, in a framework that could prove its own consistency. G\"odel showed both goals unattainable: no computable axiomatic framework proves every true statement about a sufficiently rich computational system, and no consistent framework of that strength proves its own consistency. Interpretations have mainly dwelt on the resulting lack of certainty, since on the standard view the truth of frameworks like Peano arithmetic is taken for granted. Ours places the significance elsewhere, on the incompleteness itself. As beings of finite computational power we are restricted to computable frameworks, and so to frameworks that leave some truths about any sufficiently rich system beyond reach. Mathematics can demonstrate its own limits in a way that physics, so far, cannot.

Not all theories suffer this limitation. Tarski showed that the theory of real closed fields, essentially the axioms of the ordered real line, is decidable: every statement in its language can be proved or disproved. The difference is expressive power. Real closed fields cannot encode arithmetic on natural numbers; they cannot distinguish integers from other reals, so they cannot describe systems rich enough for self-reference. Incompleteness is not a universal feature of axiomatic frameworks but a consequence of computational richness.

We can prove such a thing about our own theories because the tools of investigation and the objects of study are the same kind of thing: both are computations. Mathematics can turn its instruments on itself. This self-reference produces the limitation, but it also makes mathematics inexhaustible. Each new framework is itself a computational system, generating new traces, new patterns, and new questions that demand yet another framework. Sufficiently rich systems are always more complex than any theory about them, and the experimental activity of running computations and observing their traces can never be fully replaced by deriving consequences from axioms. This is also the moral Chaitin draws from his information-theoretic reading of G\"odel \cite{chaitin2005metamath}.

What, then, is the status of the consistency of Peano arithmetic? We can have no absolute certainty about it, yet it commands near-universal confidence. Through the eyes of this essay its status is exactly that of a law of nature. We believe no computation will ever derive a contradiction from the Peano axioms the way we believe unsupported objects fall: not because we possess a guarantee but because the traces overwhelmingly support it. At the root the evidence is enumerative induction: an enormous body of derivations run without ever surfacing a contradiction.

There is a second route to that confidence, and the axiom of infinity exemplifies it. Peano arithmetic cannot prove its own consistency, but a different system, ZF, can prove the consistency of PA. ZF encodes a different set of intuitions, about sets and the infinite, and is internally consistent as far as mathematicians have been able to test. By positing the totality of the natural numbers as a single completed infinity, ZF can check that Peano's axioms hold of that object: this is the sense in which ZF explains PA, and it is exactly the move of positing a field to explain Faraday's observations. The explanation is not an empty formality; it commits ZF to something, since a contradiction in arithmetic would refute set theory too. The posit also pays off in the purely finitary world. The known proofs of many statements about finite numbers, Fermat's Last Theorem and the Weil conjectures in full generality among them, run through sheaf cohomology and analysis, machinery far above arithmetic. Even if these theorems could in principle be proved within PA, the language of PA is not the explanatory one, and the proofs were found through intuition built by playing with the infinite. For some statements no proof within PA exists at all: Paris and Harrington showed that a finite variant of Ramsey's theorem, about colorings of finite sets, cannot be proved in PA but follows in a few lines from the infinite Ramsey theorem.

Our trust in PA is thus supported by our trust in ZF, which rests in turn on direct experience with ZF and on its success in explaining arithmetic: a virtuous circle, not the clean tower of proof and truth mathematics is often imagined to be.

Gentzen's reduction of the consistency of arithmetic to induction along the ordinal $\varepsilon_0$ is a route of the same kind, seemingly more powerful since it connects the question to a framework with its own accumulated evidence. On reflection, though, such reductions do not create evidence so much as transport it. The network of interpretations between theories is a system of conduits through which trace-evidence flows from well-tested systems to new ones, and the currency at the root is always computational experience.

This view also explains the graded confidence mathematicians actually hold. Essentially no one doubts the consistency of Peano arithmetic; the consistency of ZFC is believed a little less firmly; that of ZFC with a supercompact cardinal less still. If these were known truths the ordering would be inexplicable; for well-supported empirical laws it is just what one expects. The ordering is one of exposure: set theory predicts that arithmetic will never yield a contradiction, a supercompact cardinal predicts the same of set theory, and each new prediction is one more place the traces could say no. Such expectations are occasionally refuted, as physical laws are overthrown: Reinhardt cardinals were a seriously proposed extension of the large cardinal hierarchy until Kunen derived a contradiction in ZFC.

If mathematics were a purely logical edifice in which each step depends on every previous one, this incompleteness would be cause for alarm: errors in published proofs would be catastrophic, and the permanent gap between theory and truth would leave the enterprise on uncertain ground. In practice neither worry materializes. Mathematics is extraordinarily resilient to errors in individual papers, because the primary artifact of mathematical work is not the formal proof but the intuitions mathematicians develop about the systems they study. These intuitions let us judge whether a claimed result is consistent with the traces we have seen and reasonable to expect. Important theorems, moreover, do not sit inert; they are used, repeatedly, as inputs to further computation, and a theorem making false claims would produce wrong predictions downstream, as a false physical law eventually yields a failed experiment. Mathematics is reliable not because it is logically perfect but because it is experimental: computing within the systems our theorems describe is a constant, distributed check that no single proof, however long, could provide.

The experimental view is meant as an account of how mathematicians work, but mathematics is also a cornerstone of science, especially physics. As Galileo wrote in \textit{The Assayer}, the grand book of the universe ``is written in the language of mathematics.'' Why should the patterns we see in computational traces describe the physical world so well?

\section*{\centering{The reasonable effectiveness of mathematics}}

Wigner posed the question pointedly \cite{wigner1990unreasonable}: structures developed for purely internal reasons keep turning out to be what physics needs. Riemann developed differential geometry to study the intrinsic properties of curved surfaces; sixty years later Einstein found it the natural language of gravitation. Cardano introduced complex numbers to solve cubics he could not handle over the reals; centuries later they proved essential to quantum mechanics. Why should this be?

The question has two parts. The first is why the world is understandable at all: even a world running definite rules need not be describable by \emph{short} computational systems, and yet ours is. Maxwell needed four equations to explain every classically known electromagnetic phenomenon, from magnets and static electricity to light itself. This appears to be an empirical observation that also holds in mathematics: Gauss's one-line law for the distribution of primes is compressibility found inside arithmetic. A partial explanation is that both fields select their objects: mathematicians choose the systems that reward compression, and the celebrated core of physics collects the phenomena that compress. Both also have large domains where regularity is hard to find, such as most questions about specific Diophantine equations, or about complex systems in physics.

The second part we can address: granting this miracle, why are the structures mathematicians develop in advance, for internal reasons, so often the ones physics requires? Take the complex numbers. Introduced to solve polynomial equations, they slowly became central because they explain many unrelated phenomena, from intersections of geometric figures to representation theory and even number theory, as in several proofs of quadratic reciprocity. We may not understand precisely \emph{why} they are so central, to the extent the question can be answered at all, but that they are central across so many computational systems is itself an experimentally observed phenomenon.

This converts the apparent miracle into the selection principle noted above. Encoding a concept in an axiomatic system frees it of its original context and gives it a life of its own; the concepts reused most across computational systems are the ones that survive, become central, and acquire the best-developed theory. The physical universe, observably rich in its traces and remarkably cheap in its rules, sits squarely in the class of systems mathematicians have always curated, so it is no surprise that the concepts physics needs are the ones this selection has already favored.

The selection principle also explains the converse: why physical intuition is so unreasonably effective in mathematics. Witten won a Fields Medal for proving and re-proving deep mathematical results from physical intuition, and he is far from alone. That quantum field theory has something to say about natural numbers stops being surprising once both fields are seen to be in the business of generality, building unified frameworks that capture many processes: intuition flows both ways because the two are, in the end, the same activity.

Wigner's puzzle, on our reading, reduces to this: why are the rules of our universe short? This shortness may not be universal even in our universe; the laws of biology seem to have a very different character, with the rules generating the process being as complicated as the patterns they produce. The shortness of nature's rules is itself a pattern in our traces, one still waiting for its framework.

Mathematics, on the view we have arrived at, is the experimental study of computational systems, and physics the same activity aimed at a different target. The objects, as Demazure reminds us, are stubborn: always more general than the frameworks we build to contain them. The product, as Thurston wrote in the lines that opened this essay, is clarity and understanding, and it lives less in the frameworks than in the intuitions that exploring them trains. Those intuitions always run ahead of our formalisms, yet they are made of nothing but experience with formalisms: each new framework, once built, is the instrument for seeing the next. The cave of our revised allegory has no outside. Each exit opens into a larger cave, and the inference from shadows to shapes begins again.

\bibliography{AMSNotices}

\end{document}